\documentclass{scrartcl}

\usepackage{amsmath}
\usepackage{amsfonts}
\usepackage{amssymb} 
\usepackage{amsthm}
\usepackage{dsfont}
\usepackage{tikz}
\usepackage[symbol]{footmisc}
\usepackage{bm}

\usetikzlibrary{decorations.pathreplacing}
\newtheorem{theorem}{Theorem}
\newtheorem{lemma}{Lemma}
\newtheorem{definition}{Definition}
\newtheorem{ass}{Assumption}
\newtheorem{remark}{Remark}
\newcommand{\rvline}{\hspace*{-\arraycolsep}\vline\hspace*{-\arraycolsep}}

\newcommand{\bk}{\mathbf{k}}

\newcommand{\balpha}{\boldsymbol{\alpha}}

\newcommand{\bx}{\mathbf{x}}

\newcommand{\bz}{\mathbf{z}}

\newcommand{\bM}{\mathbf{M}}
\newcommand{\bL}{\mathbf{L}}

\newcommand{\bX}{\mathbf{X}}

\newcommand{\D}{{\mathcal{D}}}

\newcommand{\F}{{\mathcal{F}}}

\newcommand{\Nu}{{\mathcal{N}}}

\newcommand{\N}{\mathbb{N}}

\newcommand{\R}{\mathbb{R}}
\newcommand{\Z}{\mathbb{Z}}
\newcommand{\Rd}{\mathbb{R}^d}
\newcommand{\IND}{\mathbb{I}}

\DeclareMathOperator*{\argmin}{arg\,min}
\DeclareMathOperator{\sgn}{sgn}

\newcommand{\EXP}{{\mathbf E}}
\newcommand{\PROB}{{\mathbf P}}

\usepackage[round]{natbib}

\begin{document}

 {\LARGE \textbf{
    Statistical theory for image classification using deep
    convolutional neural network with cross-entropy loss
    under the hierarchical max-pooling model
    }\footnote{
Running title: {\textit{Statistical theory for image classification}}
  }}
\vspace{0.5cm}
\begin{center}
Michael Kohler$^1$ and Sophie Langer$^{2,}$\\
\textit{ 
$^1$Fachbereich Mathematik, Technische Universit\"at Darmstadt,
Schlossgartenstr. 7, 64289 Darmstadt, Germany,
email: kohler@mathematik.tu-darmstadt.de}
\\
\textit{ 
$^2$Faculty of Electrical Engineering, Mathematics and Computer Science,  University of Twente,
Drienerlolaan 5, 7522 NB Enschede, The Netherlands\\
email: s.langer@utwente.nl}
\end{center}

\begin{abstract}
Convolutional neural networks (CNNs) trained with cross-entropy loss
 have proven to be extremely successful in classifying images. 
In recent years, much work has been done to also improve 
the theoretical understanding of neural networks.  Nevertheless, it seems limited 
when these networks are trained with cross-entropy loss, mainly
because of the unboundedness of the target function.  In this paper, we aim to
fill this gap by analyzing the rate of the excess risk of a CNN classifier
trained by cross-entropy loss. Under suitable assumptions on the smoothness and 
    structure of the a posteriori probability, it is shown
    that these classifiers achieve a rate of convergence
    which is independent of the dimension of the image. 
 These rates are in line with the practical observations about CNNs.
\end{abstract}

\section{Introduction}
Deep convolutional neural networks (CNNs) have led to state-of-the-art performance in solving various problems, 
especially visual recognition tasks, see, e.g., \cite{LBH15}, \cite{K12}, \cite{Sch15},
\cite{RZ17}.  While deep learning applications are characterised above all by a high degree of flexibility, ranging from different initialisation strategies (\cite{JMLR:v10:larochelle09a}) to the choice of the right activation function (\cite{Janocha2017OnLF}) and the application of a proper learning algorithm (\cite{10.5555/3104482.3104516}), one thing has so far been chosen as fixed for classification: The \textit{log} or \textit{cross entropy} loss.  The smoothness of this loss function simplifies the optimisation procedure and shows good practical performance (\cite{GBC16}, \cite{Simonyan15}).  However,  statistical risk bounds for neural networks trained with logistic loss only exist for very restrictive conditions (see, e.g., \cite{KOK19}), mainly because of the unboundedness of the corresponding excess $\varphi$-risk minimizer (see \eqref{minimizer}),  which leads to slow convergence rates. 
\\
In general, many results on CNNs are based on considering them as a special type of feedfoward neural networks (FNNs) and then using results on FNNs to derive theoretical properties (\cite{oono2019approximation} and the literature cited therein).  Unfortunately, these results do not demonstrate situations, where CNNs outperform FNNs, which is the case in many practical applications, especially in image classification.  Generalisation bounds for CNNs with arbitrarily ordered fully connected and convolutional layers were derived in \cite{https://doi.org/10.48550/arxiv.1910.01487}.  Here the model complexity is bounded by the norm of the convolutional weights leading to tighter bounds than existing bounds for FNNs.  \cite{Yarotsky_2021} obtained approximation properties of deep CNNs, but only in an abstract setting, where it is unclear how to apply those results. \cite{KKW22} analysed plug-in classifiers based on CNNs and showed that under proper assumptions on the structure of the a posteriori probability, suitable defined CNNs trained by squared loss
achieve a rate of convergence which does not depend on the input dimension of the image.  But as, e.g., experimental results in \cite{GDN13} show, CNNs learned by cross entropy loss allow to find a better local optimum than the squared loss criterion. Therefore, CNNs learned by cross entropy loss are of higher practical relevance. 
 \\
Cross-entropy loss or, more generally, convex surrogate loss functions have been studied in \cite{B06} and \cite{10.2307/3448494}. \cite{B06} showed that for convex loss functions $\phi$ satisfying a certain uniform strict convexity condition, the rate of convergence can be strictly faster than the classical $n^{-1/2}$ , depending on the strictness of convexity of $\phi$ and the complexity of the class of classifiers. \cite{10.2307/3448494} analysed how close the optimal Bayes error rate can be approximately reached using a classification algorithm that computes a classifier by minimizing a convex upper bound of the classification error function. Some results of this article (see Lemma 1) are also used in our analysis. In Lemma 1 b) we derive a modification of Zhang's bound
which enables us to derive better rate of convergence
under proper assumptions on the a posteriori probability. 
\\
\\
In this paper we derive dimension-free rates for CNN classifiers with cross-entropy loss in a binary image classification problem, where we impose some hierarchical structure on the a posteriori probability.  In case that with high probability the a posteriori probability is very close to zero or one,  meaning that the optimal classification rule makes only a very small error, our rate can even been improved. The first result can be framed as an extension of the analysis of \cite{KKW22}, which analysed plug-in classifiers based on a class of CNNs in a similar setting.  However, it is not straightforward to extend these results to CNNs with cross-entropy loss as one cannot analyse the classification problem as a nonparametric regression setting and a network approximation for the logistic function is needed.  We deal with these difficulties with novel approximation results as well as an alternative proof strategy to bound the excess risk of the classifier.  

\section{Problem Setting}
\subsection{Image classification}
We consider a binary classification problem.  Let $d_1,d_2 \in \N$, $\mathcal{X} =[0,1]^{d_1 \times d_2}$ be an image space and $\mathcal{Y}=\{-1,1\}$ the set of corresponding binary labels. We describe an (random)
image from (random) class $Y \in \mathcal{Y}$ by a (random) matrix $\mathbf{X} \in \mathcal{X}$
with $d_1$ columns and $d_2$ rows, which contains at position $(i,j)$
the grey scale value of the pixel of the
image at the corresponding position.  Let $\PROB$ be the probability measure of $\mathcal{X} \times \mathcal{Y}$ and define by
\begin{align*}
\eta(\bx) = \PROB(Y=1|\bX=\bx)
\end{align*}
the so--called \textit{a posteriori probability}. 
\\
Our aim is to predict $Y$ by a deterministic function $g: \mathcal{X} \to \R$ such that
the sign of $g(\bX)$ is a \textit{good} prediction of $Y$.  In particular, we aim to minimize the \textit{prediction error} or $0$-$1$ \textit{risk}
\begin{align*}
\mathcal{E}(g) = \PROB(Y\sgn(g(\mathbf{X})) \leq 0) = \EXP(\mathds{1}(\sgn(g(\bX)) \neq Y))
\end{align*}
where $\sgn(x) = 1$ if $x>0$ and $-1$ otherwise and $\mathds{1}(E)$ is the the indicator function of the set $E$, that is, $1$ if event $E$ occurs and $0$ otherwise. It is well-known, that the Bayes classifier $f^*(\bx) = 2\eta(\bx) -1$ minimizes $\mathcal{E}$ among all measurable functions (cf., e.g., Theorem 2.1 in \cite{DGL96}).
But, as the probability measure $\PROB$ of $(\bX,Y)$ is unknown in practice,
we cannot find $f^*$.  Instead we estimate
$f^*$ by using the training data $\D_n = \{(\bX_i, Y_i)\}_{i=1}^n$, 
where $(\bX_i, Y_i)$ are independent copies of the random vector $(\bX, Y) \sim \PROB$.
A popular approach is estimating $f^*$ by minimizing the empirical risk
\begin{align*}
\mathcal{E}_n(g) = \frac{1}{n} \sum_{i=1}^n \mathds{1}\{\sgn(g(\bX_i)) \neq Y_i\}
\end{align*}
among a class of real-valued functions $\mathcal{F}_n$.  However
 minimizing the empirical risk with $0$-$1$ loss over $\mathcal{F}_n$ is NP hard and thus computationally not feasible (\cite{B06}).  By replacing the number of misclassifications by a convex surrogate loss $\varphi$, one can overcome computational problems. 
For a given loss
$\varphi$ we are searching for an estimate $\hat{f}_n \in \mathcal{F}_n$ that minimizes
\begin{align*}
\mathcal{E}^{\varphi}_n(g) = \frac{1}{n} \sum_{i=1}^n \varphi\left(Y_ig(\bX_i)\right).
\end{align*}
By the law of large numbers,  the empirical $\varphi$-risk converges to the population $\varphi$-risk
\begin{align*}
\mathcal{E}^{\varphi}(g)=\EXP(\varphi(Yg(\bX)))
\end{align*}
when $n \to \infty$.
A wide variety of classification methods are based on the idea of replacing the 
$0$-$1$ risk by some kind of convex surrogate loss. 
In particular, AdaBoost (\cite{Friedman_2000}) employs the \textit{exponential} loss $\exp(-z)$, while support vector machines often use  \textit{hinge} loss $\max(1-z,0)$ (\cite{Vapnik1998}) and logistic regression applies the \textit{log} loss $ \log(1+\exp(-x))$ (\cite{Hastie_2009}). 
In the context of CNNs and image classification it is a standard to use \textit{cross-entropy} loss or \textit{log} loss. 
Therefore we fix $\varphi(x) =  \log(1+\exp(-x))$ in the following.  \\
The classification performance of an estimator 
\begin{align}
\label{ERM}
\hat{f}_n \in \argmin_{g \in \mathcal{F}_n} \mathcal{E}^{\varphi}_n(g)
\end{align}
 is measured by its \textit{excess risk}
\begin{align*}
\mathcal{E}(\hat{f}_n, f^*)=\mathcal{E}(\hat{f}_n) - \mathcal{E}(f^*). 
\end{align*}
Accordingly we denote the \textit{excess $\varphi$-risk} by
\begin{align*}
\mathcal{E}^{\varphi}(\hat{f}_n, f_{\varphi}^{*})=\mathcal{E}^{\varphi}(\hat{f}_n) - \mathcal{E}^{\varphi}(f_{\varphi}^{*}),
\end{align*}
where 
\begin{align}
\label{minimizer}
f_{\varphi}^{*} = \argmin_{g \in \mathcal{F}_n} \mathcal{E}^{\varphi}(g) = \log\left(\frac{\eta(\bx)}{1-\eta(\bx)}\right)
\end{align}
in case of $\varphi(x) =  \log(1+\exp(-x))$ (cf., \cite{Friedman_2000}). 
Our following result states a relation between the excess risk and its logistic surrogate counterpart.
\begin{lemma}
  \label{le1}
  Define $\hat{f}_n$, $f^*$ and $f_{\varphi}^{*}$ as above.\\
\textbf{a)} Then
\begin{eqnarray*}
&&
\EXP\{\mathcal{E}(\hat{f}_n, f^*)\} \leq \sqrt{2} \cdot \EXP\{\mathcal{E}^{\varphi}(\hat{f}_n, f^*_\varphi)^{1/2}\}
\end{eqnarray*}
holds.  \\
\textbf{b)}
Then
\begin{eqnarray*}
&&
\EXP\{\mathcal{E}(\hat{f}_n, f^*)\} \leq 2 \cdot \EXP\{\mathcal{E}^{\varphi}(\hat{f}_n, f^{*}_{\varphi})\} + 4 \cdot \mathcal{E}^{\varphi}(f^*_{\varphi})
\end{eqnarray*}
holds. 
\\
In both parts the expectation is taken over the training data $\mathcal{D}_n$.
  \end{lemma}
\begin{remark}
We use two different bounds on the excess risk as in case of proper assumption on the distribution of $(\mathbf{X},Y)$ (see Assumption 2 below), $\mathcal{E}^{\varphi}(f^*_{\varphi})$ is small, such that part (b) of this lemma leads to faster rates.
\end{remark}
Part (a) follows from Theorem 2.1 in \cite{10.2307/3448494}, where we choose $s=2$ and $c=2^{-1/2}$.
For part (b) we set $\bar{f}_n(\bx) := 1/(1+\exp(-\hat{f}_n(\bx))$ and $g(z) := \log(z/(1-z))$ for $z \in (0,1)$.  One can show that
\begin{align}
\label{le1eq1}
\EXP\{\mathcal{E}(\hat{f}_n, f^*)\} &\leq 2 \cdot \EXP\{|\bar{f}_n(\bX) - \eta(\bX)|\}.
\end{align}
For $h_1(z) := \varphi(1 \cdot g(z))$ and $h_2(z):= \varphi(-1 \cdot g(z))$ it further holds that
$|h'_j(z)| \geq 1$ for $j \in \{1, 2\}$ and $z \in (0,1)$.
Using mean value theorem we can bound \eqref{le1eq1} by
\begin{align*}
&2 \cdot \EXP\bigg\{\left|\varphi\left(Y \cdot g(\bar{f}_n(\bX))\right)- \varphi\left(Y \cdot g(\eta(\bX))\right)\right|\bigg\}.
 \end{align*}
 With $|a-b| \leq a+b$ for $a,b \geq 0$, the assertion follows. The complete proof is found in the appendix.

\subsection{Hierarchical max-pooling model}
In order to derive nontrivial rate of convergence results
on the excess $\varphi$-risk of any estimate it is necessary
to restrict the class of distributions (cf.,
\cite{C68} and \cite{D82}).  In case of logistic loss we have 
$f_{\varphi}^{*}(\bx)  = \log(\eta(\bx)/(1-\eta(\bx)))$, showing that $f_{\varphi}^{*}$
is a monotone transformation of the a posteriori probability $\eta$.  Hence we need
to impose some assumptions on $\eta$. 
\\
As in \cite{KKW22} we assume that
the a posteriori probability fulfills some $(p,C)$\textit{-smooth hierarchical max-pooling model}. 
As smoothness measure we use the following definition of  $(p,C)$-smoothness. 
For simplicity we introduce the multi-index notation, that is,  $\partial^{\balpha} = \partial^{\alpha_1} \dots \partial^{\alpha_d}$ with
$\balpha = (\alpha_1, \dots, \alpha_d) \in \N_0^d$.

\begin{definition}
\label{de3}
Let $p=q+s$ for some $q \in \N_0$ and $0< s \leq 1$.
A function $f:\R^d \rightarrow \R$ is called
\textbf{$(p,C)$-smooth}, if for every $\balpha \in
\N_0^d$
with $\|\balpha\|_1 = q$ the partial derivative $\partial^{\balpha} f$
exists and satisfies
\[
|\partial^{\balpha} f(\bx) - \partial^{\balpha} f(\mathbf{z})| \leq C \cdot \|\bx - \bz\|^s
\]
for all $\bx,\bz \in \R^d$.
\end{definition}

For the next definitions we frequently use the following notation:
For $M \subseteq \R^d$ and $\bx \in \R^d$ we define
\[
\bx+M = \{\bx+\bz \, : \, \bz \in M \}.
\]
For $I \subseteq d_1 \times d_2$
and
$\bx=(x_i)_{ i \in d_1  \times d_2}
\in  [0,1]^{d_1 \times d_2}$ we set
\[
\bx_I =(x_i)_{i \in I}.
\]
  
  The definition of \textit{hierarchical max-pooling models} is motivated by the following observation: Human beings often decide, whether a given image contains some object, i.e., a car, or not by scanning subparts of the image and checking, whether the searched object is on this subpart. For each subpart the human estimates a probability that the searched object is on it. The probability that the whole image contains the object is then simply the maximum of the probabilities for each subpart of the image. This leads to the definition of a max-pooling model for the a posteriori probability. 
  
  \begin{definition}
  \label{de1}
	Let $d_1,d_2\in\N$ with $d_1,d_2>1$ and $m: [0,1]^{d_1\times d_2} \rightarrow \R$.
We say that $m$
satisfies a \textbf{max-pooling model with index set}
\[
I \subseteq d_1-1 \times d_2-1,
\]
if there exist a function $f:[0,1]^{(1,1)+I}\rightarrow \R$ such that
\[
m(\bx)=
\max_{
  (i,j) \in \Z^2 \, : \,
  (i,j)+I \subseteq d_1 \times d_2
}
f\left(
\bx_{(i,j)+I}
\right)
\]
for $\bx \in [0,1]^{d_1 \times d_2}$.
\end{definition}
Additionally, the probability that a subpart contains the searched object is composed by several decisions, if parts of the searched objects are identifiable. This motivates the hierarchical structure of our model.  In the following we denote the four block matrices of a matrix $\mathbf{x} \in [0,1]^{2^k \times 2^k}$ by $\mathbf{x}_{1,1}, \mathbf{x}_{2,1}, \mathbf{x}_{1,2}, \mathbf{x}_{2,2}$, where 
\begin{align*}
\mathbf{x}_{i,j} = \mathbf{x}_{\{(i-1)2^{k-1}+1, \dots,  i2^{k-1}\} \times \{(j-1)2^{k-1}+1, \dots,  j2^{k-1}\}} \in [0,1]^{2^{k-1} \times 2^{k-1}},
\end{align*}
$i,j \in \{1,2\}$. This means that
 \begin{align}
 \label{defx}
 \mathbf{x}=
\begin{pmatrix}
\mathbf{x}_{1,1} & \rvline & \mathbf{x}_{2,1}\\
\hline
\mathbf{x}_{1,2} & \rvline &  \mathbf{x}_{2,2}
\end{pmatrix}.
 \end{align}

\begin{figure}
\centering
\def\svgwidth{200pt}    
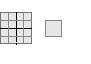  
    \caption{Illustration of Definition \ref{de2} with $l=2$}
    \label{fig:def3}
\end{figure}

\begin{definition}
\label{de2}
Let $d_1,d_2\in\N$ with $d_1,d_2>1$ and $m: [0,1]^{d_1\times d_2} \rightarrow \R$. 
    We say that
\[
f:[0,1]^{2^l \times 2^l} \rightarrow \R
\]
 satisfies a
    \textbf{hierarchical model of level $l$},
    if there exist functions
    \[
    g_{k,s}: \R^4 \rightarrow [0,1]
    \quad (k=1, \dots, l, s=1, \dots, 4^{l-k} )
    \]
    such that we have
    \[
f=f_{l,1}
    \]
    for some
    $f_{k,s} :[0,1]^{2^k \times 2^k} \rightarrow \R$ recursively defined by
    \begin{align*}
    f_{k,s}(\bx)=&g_{k,s} \big(
    f_{k-1,4 \cdot (s-1)+1}(\bx_{1,1})
    , 
        f_{k-1,4 \cdot (s-1)+2}(\bx_{2,1}), \
        f_{k-1,4 \cdot (s-1)+3}(\bx_{1,2}
        ), 
        f_{k-1,4 \cdot s}(\bx_{2,2})
    \big)
    \end{align*}
    for $k=2, \dots, l, s=1, \dots,4^{l-k}$ and $\bx \in
[0,1]^{2^k \times 2^k}$
    and
    \[
 f_{1,s}(
x_{1,1},x_{1,2},x_{2,1},x_{2,2}
)= g_{1,s}(x_{1,1},x_{1,2},x_{2,1},x_{2,2})
 \]
 for $s=1, \dots, 4^{l-1}$ and $x_{1,1},x_{1,2},x_{2,1},x_{2,2} \in [0,1]$.  
\end{definition}
An illustration of Definition \ref{de2} for $l=2$ is shown in Figure \ref{fig:def3}.

Combining Definition \ref{de3}, \ref{de1} and \ref{de2} leads to the final definition of $(p,C)$-\textit{smooth hierarchical max-pooling models}.

\begin{definition}
  \label{de4}
     We say that
     $m: [0,1]^{d_1 \times d_2} \rightarrow \R$
     satisfies a \textbf{$(p,C)$-smooth hierarchical max-pooling model of level $l$
      } (where $2^l \leq \min\{ d_1,d_2\}$),
     if $m$ satisfies a max-pooling model with index set $I=2^{l}-1 \times 2^{l}-1$,
      the function
     $f:[0,1]^{(1,1)+I} \rightarrow \R$ in the definition of this
     max-pooling model satisfies a hierarchical model
     with level ${l}$ and if all functions $g_{k,s}$ in the definition of the functions $m$ are $(p,C)$--smooth for some $C>0$.
  \end{definition}

\subsection{Convolutional neural networks}
We consider CNNs that take $d_1 \times d_2$-dimensional images as input and produce an one-dimensional output. As the name suggests, the most important operation of a CNN is its convolution.  The main idea behind it is to apply filters, i.e., small weight matrices to the input image to extract high-level information.  Mathematically a convolution can be described as follows: Let $\bX$ be a $d_1 \times d_2$ input matrix, $\bX_{i,j}$ be its $\ell \times \ell$ block matrix with entries $(X_{i+a,j+b})_{a,b = 0, \dots, \ell-1}$ and $\mathbf{W}$ be a corresponding filter of size $\ell$.  The entry $(i,j)$ of a resulting channel $\mathbf{\tilde{C}}$ can then be described by
\begin{align}
\label{tildeC}
\tilde{C}_{i,j} = \sum_{k,m=1}^{\ell} (\bX_{i,j} \odot \mathbf{W})_{k,m},
\end{align}
where $\odot$ denotes the Hadamard product.  Finally an activation function $\sigma$ is applied componentwise,  i.e., $C_{i,j} := \sigma(\tilde{C}_{i,j})$. This in turn means that the final channel $\mathbf{C}$ consists of entries computed by the sum of a Hadamard product between the filter and the respective block matrix of the input applied to an activation function $\sigma$.  We set
\begin{align*}
\mathbf{C} := \sigma(\mathbf{W} \star \mathbf{X})
\end{align*}
with $\sigma(x) = \max\{x,0\}$ being the ReLU activation function.
Here $\star$ describes the computation of each entry $\tilde{C}_{i,j}$ as in \eqref{tildeC}, where $\sigma$ is applied componentwise.  One can see that the weights generating the feature map $\mathbf{C}$ are shared, which has the advantage of reducing the complexity of the model and the training time of the networks.  Usually a CNN consists of several convolutional layers.  Each convolutional layer $l$ $(l \in \{1, \dots, L\})$ consists of $k_l \in \N$ channels (also called feature maps) while the filter size $M_l \in \{1, \dots, \min\{d_1, d_2\}\}$ per layer is fixed.  In our setting we make use of so-called zero-padding meaning that we enlarge each channel by appending zero matrices on each side such that the convolution does not change the in-plane dimension.  This, in turn, means that every resulting channel has size $d_1 \times d_2$.  For filters
\begin{align*}
(\mathbf{W}_{s_1,s_2, l})_{l = 1,\dots, L, s_1= 1, \dots, k_{l-1}, s_2 = 1, \dots, k_l}
\end{align*} 
the $s$-th channel of layer $l$ $(s=1, \dots, k_l, l=1, \dots, L)$ can be described by
\begin{align}
\label{CNN}
\mathbf{C}_{s,l} = \sigma\left(\sum_{s_1=1}^{k_{l-1}} \mathbf{W}_{s_1, s, l} \star \mathbf{C}_{s_1, l-1}\right)
\end{align}
with $\mathbf{C}_{1,0} = \bX$ and $k_0 =1$.

In our network, only in the last step a max-pooling layer is applied to the values of the last convolutional layer $L$.  As in \cite{LS22} we consider a global max-pooling where we extract from every channel $\mathbf{C}_{s,L}$ $(s=1, \dots, k_L)$ the largest absolute value.  A CNN with $L \in \N$ convolutional layers and one pooling layer, channel vector $\mathbf{k}=(k_1, \dots, k_L) \in \N^{L}$ and filter vector $\mathbf{M}=(M_1, \dots, M_L) \in \N^L$,  where $k_i$ describes the number of channels and $M_i$ describes the size of the filters in layer $i$, respectively, can be described as a funcion $f:[0,1]^{d_1 \times d_2} \to \mathbb{R}^{k_L}$ with
\begin{align*}
\mathbf{x} \mapsto f(\mathbf{x}) = (|\mathbf{C}_{1,L}|_{\infty}, \dots,  |\mathbf{C}_{k_L,L}|_{\infty})
\end{align*}
with $\mathbf{C}_{s,L}$ recursively defined as in \eqref{CNN}. We denote this network class by $\mathcal{F}_{L, \mathbf{k}, \mathbf{M}}^{C}$. 

After convolutional and pooling layers typically several fully connected layers are applied.  Again we choose the ReLU activation function $\sigma(x)=\max\{x,0\}$.  Following the definition in \cite{Sch17}, a fully connected network with $L \in \N$ hidden layers and width vector $\mathbf{k}=(k_0, \dots, k_{L+1}) \in \N^{L+2}$, where $k_i$ denotes the number of neurons in layer $i$, can be described by a function $f:\R^{k_0} \to \R^{k_{L+1}}$ with
\begin{align*}
\bx \mapsto f(\bx) = \mathbf{W}_L \sigma_{v_L} \mathbf{W}_{L-1}\sigma_{v_{L-1}} \cdots \mathbf{W}_1 \sigma_{\mathbf{v}_1} \mathbf{W}_0 \mathbf{x},
\end{align*}
where $\mathbf{W}_j$ is a $k_j \times k_{j+1}$ weight matrix and $\mathbf{v}_j$ is the $j$-th shift vector.  We denote the network class of fully connected neural networks by $\mathcal{F}_{L, \mathbf{k}}$.

Our final function class $\mathcal{F}_{n, \Theta}$ is then a composition of convolutional and fully connected layers, i.e., 
\begin{align}
\label{CNNclass}
  \mathcal{F}_{n, \Theta} = &\bigg\{g \circ f: f \in \mathcal{F}^{C}_{L_n^{(1)}, \mathbf{k}^{(1)}, \mathbf{M}}, g \in \mathcal{F}_{L_n^{(2)}, \mathbf{k}^{(2)}} 
  , \|g \circ f \|_\infty \leq
  \beta_n \bigg\}, 
\end{align}
where $\Theta = (\bL, \bk^{(1)}, \bk^{(2)}, \bM)$ with parameters
\begin{align*}&\bL=(L_n^{(1)},L_n^{(2)}), ~\bk^{(1)}=\left(k_1^{(1)},\dots,k_{L_n^{(1)}}^{(1)}\right),\\
&~ \bk^{(2)}=\left(k_1^{(2)},\dots,k_{L_n^{(2)}}^{(2)}\right),
~ \bM=(M_1,\dots,M_{L_n^{(1)}})
\end{align*}
and $\beta_n = c_1 \cdot \log n$ for some constant $c_1 >0$. Accordingly we denote by 
\begin{align}
\label{est}
\hat{f}^{CNN}_n = \argmin_{f \in \mathcal{F}_{n, \Theta}} \frac{1}{n} \sum_{i=1}^n \log(1+\exp(-Y_i \cdot f(\mathbf{X}_i)))
\end{align}
the empricial risk minimizer based on our class of CNNs.
\\
\\
In general, deep learning theory can be roughly divided into three parts, namely expressivity, generalisation and optimisation (see \cite{K20}). While the intersection of all three aspects has only been analysed in very limited settings so far, e.g. , for shallow neural networks and a restricted class of regression functions (see, e.g. , \cite{BKWL21}), most statistical risk bounds of neural networks exclude the optimisation algorithm and deal with the empirical risk minimizer (ERM) instead (see, e.g., \cite{Sch17}, \cite{BK17}, \cite{KL20}). Following this line of work, we also analyse the ERM based on a particular class of CNNs.  It therefore remains an open question whether similar rates can be shown for CNNs trained with (stochastic) gradient descent.  In case of overparametrized CNNs, e.g. , \cite{CS19} could show that the gradient descent is able to find the global minimum of the empirical loss function.  But for overparametrized networks one cannot use standard generalisation bounds as these usually depend on the number of parameters. Therefore a completely new statistical approach is needed for this analysis.

\section{Main Result}
In this section we derive convergence rates of the excess risk of $\hat{f}^{CNN}_n$ under the assumption that the a posteriori probability $\eta$ fulfills a $(p,C)$-smooth hierarchical max-pooling model (see Definition \ref{de4}). Before providing our results,  two assumptions are imposed on the distribution of $(\bX, Y)$. 
\begin{ass}
\label{ass1}
For $p \geq 1$ and $C>0$ arbitrary, 
$\eta(\bx) = \PROB\{Y=1|\mathbf{X} = \bx\}$ satisfies a $(p,C)$--smooth hierarchical
max-pooling model of finite level $l$ and $supp(\PROB_{\bX}) \subseteq [0,1]^{d_1 \times d_2}$.
\end{ass}
The second is a margin condition on the a posteriori probability.
\begin{ass}
\label{ass2}
For $f_{\varphi}^*$ being the minimizer of $\mathcal{E}^{\varphi}(g)$ and $n \in \N$, it holds
   \begin{equation*}
\PROB \left\{\bX : |f_{\varphi}^*(\bX)| > \frac{1}{2} \cdot \log n\right\} \geq
1 - \frac{1}{\sqrt{n}}.
\end{equation*}
\end{ass}

Assumption \ref{ass2} requires that with high probability
            the a posteriori probability is very close to zero or one, 
            and hence the optimal classification
            rule makes only a very small error. This is in particular realistic
            for various image classification tasks, where object classes can often be confidently distinguished (cf., \cite{KOK19}).  
            Additionally, a similar assumption was applied in \cite{BSch22} within the context of a multiclass classification problem, also analysing a ReLU network classifier minimising
            cross entropy loss.

\begin{theorem}
	  \label{th1}
	  Suppose Assumption \ref{ass1} holds.
Set
\[
L_n^{(1)}= \frac{4^l-1}{3} \cdot \lceil c_3 \cdot n^{2/(2p+4)} \rceil+l
\quad \mbox{and} \quad
L_n^{(2)}= \lceil c_2 \cdot n^{1/4} \rceil,
\]
\[
M_s = 2^{\pi(s)} \quad (s=1, \dots, L_n^{(1)}),
\]
where the function
$\pi:\{1, \dots, L_n^{(1)} \}
\rightarrow \{1, \dots, l\}$ is defined by
\[
\pi(s)=\sum_{i=1}^l
\mathds{1}_{
  \{
s \geq i+\sum_{r=l-i+1}^{l-1} 4^r \cdot \lceil c_3 \cdot n^{2/(2p+4)} \rceil
  \}
},
\]
choose $\mathbf{k}^{(1)} =(c_4, \dots, c_4) \in \N^{L_n^{(1)}}$ and
$\mathbf{k}^{(2)}=(c_5, \dots, c_5) \in \N^{L_n^{(2)}}$,
and define the estimate $\hat{f}^{CNN}_n$
as in \eqref{est}.
Assume that the constants $c_2, \dots, c_5$ are sufficiently large.

\noindent
\textbf{a)} There exists a constant $c_6 =c_6(\eta, d_1, d_2, p,C, l) >0$ such that we have for any $n > 1$
\begin{eqnarray*}
&&
\EXP\{\mathcal{E}(\hat{f}^{CNN}_n, f^*)\}
\leq
c_6 \cdot (\log n) \cdot n^{-\min\{\frac{p}{4p+8},
   \frac{1}{8} \}}.
\end{eqnarray*}
\noindent\\
    \textbf{b)}
    If, in addition, Assumption \ref{ass2} holds
, then there exists a constant $c_7 = c_7(\eta, d_1, d_2, p,C, l) >0$ such that we have for any
    $n > 1$
\begin{eqnarray*}
&&
\EXP\{\mathcal{E}(\hat{f}^{CNN}_n, f^*)\}
\leq
c_7 \cdot (\log n)^{2} \cdot n^{-\min\{\frac{p}{2p+4},
   \frac{1}{4} \}}.
\end{eqnarray*}    
In both parts the expectation is taken over the training data $\mathcal{D}_n$.
\end{theorem}

\begin{remark}   
An interesting feature of the convergence rates in Theorem \ref{th1} 
is that both rates do not depend on the dimension $d_1 \cdot d_2$ of the input image.
Thus, given the structure of the a posteriori probability fulfills a $(p,C)$-smooth 
hierarchical max-pooling model,  our estimator circumvents the curse of dimensionality.
Under Assumption \ref{ass2} the rate can even be improved.  To us these results partly 
explain the good performance of CNN classifiers on image data.
\end{remark}

\begin{remark}
            The definition of the parameters $L_n^{(1)}$ and $M_i$ $(i=1, \dots, L_n^{(1)})$
            of the estimate in Theorem \ref{th1} depends on the smoothness
            and the level of the hierarchical max-pooling model
            for the a posteriori probability, which are
            usually unknown in applications. In this case it is
            possible to define these parameters in a data-dependent
            way, e.g., by using a splitting of the sample approach
            (cf., e.g., Chapter 7 in Gy\"orif et al. (2002)).
\end{remark}
\textbf{On the proof.} 
To prove Theorem \ref{th1} we use Lemma \ref{le1} in combination with the following general upper bound on the excess $\varphi$-risk of an empirical risk minimizer $\hat{f}_n \in \mathcal{F}_n$, where $\mathcal{F}_n$ can be a general function space consisting of functions $f:\R^{d_1 \times d_2} \to \R$. 
\begin{lemma}
  \label{le2}
  Let $\varphi$ be the logistic loss and $\mathcal{D}_n=\{(\bX_i, Y_i)\}_{i=1}^n$. Then the the empirical risk minimizer $\hat{f}_n$ defined as in \eqref{ERM} satisfies 
\begin{align*}
\EXP\{\mathcal{E}^{\varphi}(\hat{f}_n,  f_{\varphi}^*)\} \leq &2
\cdot
\sup_{f \in \F_n}|\mathcal{E}^{\varphi}(f) - \mathcal{E}^{\varphi}_n(f)| + \inf_{f \in \mathcal{F}_n} \mathcal{E}^{\varphi}(f, f_{\varphi}^*).
\end{align*}
  \end{lemma}
  
Lemma \ref{le2} shows that the excess $\varphi$-risk of an ERM is bounded above by the sum of two terms.  The first term is the so-called \textit{generalisation} error. It is closely related to the complexity of the function class and can be bounded using results from empirical process theory.  The second one is the \textit{approximation error} measuring how rich the function class $\mathcal{F}_n$ is, meaning if we can express the problem under consideration by a function of $\mathcal{F}_n$.  \\
In the following $\mathcal{F}_n$ is chosen to be a class of convolutional neural networks, i.e., $\mathcal{F}_n = \mathcal{F}_{n, \Theta}$ and the estimator under consideration is defined as in \eqref{est}.

\section{Approximation error}
To bound $\inf_{f \in \mathcal{F}_{n, \Theta}} \mathcal{E}^{\varphi}(f, f_{\varphi}^*)$ we use that for an arbitrary $h \in \mathcal{F}_{n, \Theta}$ 
\allowdisplaybreaks
\begin{align*}
\inf_{f \in \mathcal{F}_{n, \Theta}} \mathcal{E}^{\varphi}(f, f_{\varphi}^*) \leq
 \mathcal{E}^{\varphi}(h, f_{\varphi}^*) 
& =
\int
(
\eta(\bx) \cdot \varphi( h(\bx))
+
(1-\eta(\bx)) \cdot \varphi(- h(\bx))
)
\PROB_\bX(d\bx)
\\
&
\quad
-
\int
(
\eta(\bx) \cdot \varphi( f_\varphi^*(\bx))
+
(1-\eta(\bx)) \cdot \varphi(- f_\varphi^*(\bx))
)
\PROB_\bX(d\bx)
\\
&
\leq
\sup_{\bx \in [0,1]^{d_1 \times d_2}}
\Big(
\left|
\eta(\bx) \cdot
\left(
\varphi(h(\bx))-
\varphi(f_\varphi^*(\bx))
\right)
\right|
\\
&
\quad
+
\left|
(1-\eta(\bx)) \cdot
\left(
\varphi(-h(\bx))-
\varphi(-f_\varphi^*(\bx))
\right)
\right|
\Big)\\
&
\leq
\sup_{\bx \in [0,1]^{d_1 \times d_2}}
\Big(
\left|
\eta(\bx) \cdot
\left(
\varphi(h(\bx))-
\varphi(g(\eta(\bx))
\right)
\right|
\\
&
\quad
+
\left|
(1-\eta(\bx)) \cdot
\left(
\varphi(-h(\bx))-
\varphi(-g(\eta(\bx)))
\right)
\right|
\Big),
\end{align*}
where
 \[
g(z)=\begin{cases}
\infty & ,z=1 \\
\log \frac{z}{1-z} & ,0<z<1 \\
- \infty & ,z=0.
\end{cases}
\]
This, in turn, means that in order to find a satisfying bound for our approximation error we need to build a CNN which approximates $g(\eta(\bx))$ properly.  Using the compositional structure of neural networks, one can break this task down into two parts. On the one hand we show that CNNs approximate $\eta(\bx)$,  i.e., $(p,C)$-smooth hierarchical max-pooling models.  On the other hand we build a fully connected neural network that approximates $g$.  The approximation result on $g$ is the following.
\begin{lemma}
  \label{le7}
  Set
  \[
g(z)=\begin{cases}
\infty & ,z=1 \\
\log \frac{z}{1-z} & ,0<z<1 \\
- \infty & ,z=0
\end{cases}
\]
and let $K \in \N$ with $K \geq 6$. Let $\eta:\Rd \rightarrow [0,1]$ and let
$\bar{\eta}:\Rd \rightarrow \R$ such that
$\|\bar{\eta}-\eta\|_\infty \leq \epsilon$ for some
$0 \leq \epsilon \leq 1/K$. Then there exists a neural
network $\bar{g}:\R \rightarrow \R$ with ReLU activation function,
$K+3$ hidden layers with $7$  neurons per layer, which is
bounded in absolute value by $\log (K+1)$ and which
satisfies
\begin{eqnarray*}
  &&
\sup_{\bx \in \R^{d_1 \times d_2}}
\left(
\left|
\eta(\bx) \cdot
\left(
\varphi(\bar{g}(\bar{\eta}(\bx))-
\varphi(g(\eta(\bx))
\right)
\right|\right.
\\
&&\quad \left.
+
\left|
(1-\eta(\bx)) \cdot
\left(
\varphi(-\bar{g}(\bar{\eta}(\bx))-
\varphi(-g(\eta(\bx))
\right)
\right|
\right)
\\
&&
\leq
c_{10} \cdot \left(\frac{\log K}{K} + \epsilon\right).
\end{eqnarray*}
\end{lemma}
A complete proof is found in the appendix. Roughly, it is based on the idea that functions of the form
\begin{align*}
\bar{g}(z) := \sum_{k=-1}^{K+1} a_k \cdot B_k(z),
\end{align*}
where
 \[
    B_k(z) =
    \begin{cases}
      0 & ,z < \frac{k-1}{K} \\
     K \cdot( z- \frac{k-1}{K})  & ,\frac{k-1}{K} \leq z < \frac{k}{K} \\
     K \cdot (\frac{k+1}{K}-z)  & ,\frac{k}{K} \leq z < \frac{k+1}{K} \\
     0 &, z \geq \frac{k+1}{K},
      \end{cases}
    \]
  can be computed by a ReLU network with $K+3$ hidden layers and 
    $7$ neurons per layer.  
    \\
    \\
Combining Lemma \ref{le7} with the approximation result on the hierarchical max-pooling models, 
 we are then able to show the following approximation result.
\begin{theorem}
  \label{le5}
  Suppose Assumption \ref{ass1} holds.
Let $\F_{n, \Theta}$ be the set of all CNNs 
with $\Theta= (\bL, \bk^{(1)}, \bk^{(2)}, \mathbf{M})$,
where $\bk^{(1)}=(c_4, \dots, c_4)$ with $c_4$ sufficiently large
and $\bk^{(2)}=(7, \dots, 7)$.

  Furthermore  assume $(L_n^{(1)})^{2p/d}
  \geq c_{8} \cdot L_n^{(2)}$.
Then
\begin{eqnarray*}
&&
\inf_{f \in \F_{n, \Theta}}
\mathcal{E}(f, f^*_{\varphi})
\leq c_{9} \cdot \left(
\frac{\log L_n^{(2)}}{L_n^{(2)}} + \frac{1}{(L_n^{(1)})^{2p/4}}
\right), 
\end{eqnarray*}
with constant $c_9=c_9(\eta, p, C,l)$.
  \end{theorem}
  The complete proof of this result is given in the appendix.

\section{Generalisation error}
The generalisation error $\sup_{f \in \F_n}|\mathcal{E}^{\varphi}(f) - \mathcal{E}^{\varphi}_n(f)|$ can be bounded using results from empirical process theory together with bounds on the covering number of CNNs.  

In particular, using Theorem 9.1 in Györfi et al. (2002) it holds for
 \begin{align*}
\mathbf{Z}=(\bX, Y), \mathbf{Z}_1=(\bX_1, Y_1), \dots, \mathbf{Z}_n=(\bX_n, Y_n), 
\end{align*}
and $\epsilon >0$, that
 \begin{align*}
\PROB\left\{\sup_{f \in \mathcal{F}_n} \left|\mathcal{E}^{\varphi}(f) - \mathcal{E}_n^{\varphi}(f)\right| > \epsilon\right\}
& = \PROB\left\{\sup_{h \in \mathcal{H}_n} \left|\EXP h(\mathbf{Z}) - \frac{1}{n} \sum_{i=1}^n h(\mathbf{Z}_i)\right| > \epsilon\right\}\\
& \leq 8 \EXP\left\{\mathcal{N}_1\left(\frac{\epsilon}{8}, \mathcal{H}_n, \mathbf{Z}_1^n\right)\right\}e^{-\frac{n \epsilon^2}{128 \cdot c_{10}^2 \cdot (\log n)^2}}.
 \end{align*}
Here 
 \begin{align*}
  \mathcal{H}_n=& \{h: \R^{d_1 \times d_2} \times \R \to \R: \exists f \in
  \mathcal{F}_{n, \Theta} \ \text{such that} \ h(\bx,y) = \varphi(y\cdot f(\bx))\}
\end{align*}
and $\Nu_1(\varepsilon,\F,\bx_1^n)$ describes the  $\varepsilon$--covering number of $\F$ on $\bx_1^n$, that is 
the smallest $\varepsilon$--cover
  of $\F$ on $\bx_1^n$, i.e., the number $N \in \N$ such that there exists  $i \in \{1, \dots, N\}$
  such that
  \[
\frac{1}{n} \sum_{k=1}^n |f(\bx_k)-f_i(\bx_k)| < \varepsilon.
  \]
  As every $\epsilon$-cover of $\mathcal{F}_{n, \Theta}$ is an $\epsilon$-cover of $\mathcal{H}_n$, we have 
\begin{align*}
\mathcal{N}_1\left(\frac{\epsilon}{8}, \mathcal{H}_n, \mathbf{Z}_1^n\right) \leq \mathcal{N}_1\left(\frac{\epsilon}{8}, \mathcal{F}_{n, \Theta}, \mathbf{X}_1^n\right).
\end{align*}
The following bound on the covering number of $\mathcal{F}_{n, \Theta}$,  then helps us to bound the generalisation error.

\begin{lemma}
  \label{le4}
Define $\mathcal{F}_{n, \Theta}$ as in \eqref{CNNclass}
and set
\begin{align*}
&k_{max}=\max\left\{k_1^{(1)}, \dots, k_{L_n^{(1)}}^{(1)},k_1^{(2)},\dots,k_{L_n^{(2)}}^{(2)}\right\},\\
&M_{max}=\max\{ M_1, \dots, M_{L_n^{(1)}} \}
\end{align*}
and
\[L_{max}=\max\{L_n^{(1)},L_n^{(2)}\}.\]
Assume $d_1\cdot d_2>1$ and $\beta_n=c_1 \cdot \log n \geq 2$.
Then we have for any
$\epsilon \in (0,1)$,
\begin{eqnarray*}
&&
\sup_{\bx_1^n \in (\R^{ d_1 \times d_2})^n} \log\left(
 \mathcal{N}_1 \left(\epsilon,\F_{n, \Theta}, \bx_1^n\right) \right)
\\
&&
\leq
c_{11} \cdot L_{max}^2 \cdot \log(L_{max} \cdot d_1 \cdot d_2) \cdot
\log \left(
\frac{c_1\cdot\log n}{\epsilon}
\right)
\end{eqnarray*}
for some constant $c_{11} >0$ which depends only on $k_{max}$ and $M_{max}$.
  \end{lemma}
The proof of this result follows by Lemma 7 in \cite{KKW22}.

\bibliography{Literatur}

\appendix
\section{FURTHER PROOFS}

\subsection{Proofs of Section 2}

In this section, we present the detailed proof of Lemma 1.
\\
\\
\noindent
 \textbf{Proof of Lemma 1.}
 \textbf{a)} This result follows from Theorem 2.1 in Zhang (2004), where we choose $s=2$ and $c=2^{-1/2}$. \\
 \textbf{b)} Set $\bar{f}_n(\bx) = 1/(1+\exp(-\hat{f}_n(\bx))$. Then we have
 \begin{align}
 \label{le2eq1}
 \EXP\{\mathcal{E}(\hat{f}_n, f^*)\} 
& \leq \EXP\Big\{\int
\Big(
(1-\eta(\bx)) \cdot \mathds{1}_{\{\hat{f}_n(\bx) \geq 0\}}(\bx) + \eta(\bx) \cdot \mathds{1}_{\{\hat{f}_n(\bx) < 0\}}(\bx) \notag\\
& \quad \quad - (1-\eta(\bx)) \cdot \mathds{1}_{\{\eta(\bx) \geq \frac{1}{2}\}}(\bx) + \eta(\bx) \cdot \mathds{1}_{\{\eta(\bx) < \frac{1}{2}\}}(\bx)
\Big)
\PROB_\bX(d\bx)\Big\}\notag\\
&  = \EXP\Big\{\int \Big(
(1-\eta(\bx)) \cdot \mathds{1}_{\{\hat{f}_n(\bx) \geq 0\}}(\bx) + \eta(\bx) \cdot \mathds{1}_{\{\hat{f}_n(\bx) < 0\}}(\bx) \notag\\
& \quad \quad - (1-\bar{f}_n(\bx)) \cdot \mathds{1}_{\{\hat{f}_n(\bx) \geq 0\}}(\bx) - \bar{f}_n(\bx) \cdot \mathds{1}_{\{\hat{f}_n(\bx) < 0\}}(\bx) \notag\\
& \quad \quad + (1-\bar{f}_n(\bx)) \cdot \mathds{1}_{\{\hat{f}_n(\bx) \geq 0\}}(\bx)
+ \bar{f}_n(\bx) \cdot \mathds{1}_{\{\hat{f}_n(\bx) < 0\}}(\bx) \notag\\
& \quad \quad -  (1-\bar{f}_n(\bx)) \cdot \mathds{1}_{\{\eta(\bx) \geq \frac{1}{2}\}}(\bx) - \bar{f}_n(\bx) \cdot \mathds{1}_{\{\eta(x) < \frac{1}{2}\}}(\bx) \notag\\
& \quad \quad+  (1-\bar{f}_n(\bx)) \cdot \mathds{1}_{\{\eta(\bx) \geq \frac{1}{2}\}}(\bx) + \bar{f}_n(\bx) \cdot \mathds{1}_{\{\eta(\bx) < \frac{1}{2}\}}(\bx) \notag\\
& \quad \quad - (1-\eta(\bx)) \cdot \mathds{1}_{\{\eta(\bx) \geq \frac{1}{2}\}}(\bx)
- \eta(x) \cdot \mathds{1}_{\{\eta(\bx) < \frac{1}{2}\}}(x)
\Big)
\PROB_\bX(d\bx)\Big\} \notag\\
& = \EXP\Big\{\int \big(
(1-\eta(\bx)-1+\bar{f}_n(\bx)) \cdot \mathds{1}_{\{\hat{f}_n(\bx) \geq 0\}}(\bx) + (\eta(\bx) - \bar{f}_n(\bx)) \cdot \mathds{1}_{\{\hat{f}_n(\bx) < 0\}}(\bx) \notag\\
& \quad \quad + (1-\bar{f}_n(\bx)-1+ \eta(\bx)) \cdot \mathds{1}_{\{\eta(\bx) \geq \frac{1}{2}\}}(\bx) + (\bar{f}_n(\bx)-\eta(\bx)) \cdot \mathds{1}_{\{\eta(\bx) < \frac{1}{2}\}}(\bx)\notag\\
& \quad \quad +  (1-\bar{f}_n(\bx)) \cdot \mathds{1}_{\{\hat{f}_n(\bx) \geq 0\}}(\bx) + \bar{f}_n(\bx) \cdot \mathds{1}_{\{\hat{f}_n(\bx) <  0\}}(\bx) \notag\\
& \quad \quad - (1-\bar{f}_n(\bx)) \cdot \mathds{1}_{\{\eta(\bx) \geq \frac{1}{2}\}}(\bx) - \bar{f}_n(\bx) \cdot \mathds{1}_{\{\eta(\bx) <  \frac{1}{2}\}}(\bx)
\big)
\PROB_\bX(d\bx)\Big\}\notag\\
& \leq 2 \cdot \EXP\{|\bar{f}_n(\bX) - \eta(\bX)|\},
 \end{align}
 where the expectation is taken over the training data $\mathcal{D}_n$. Note that the last inequality follows since 
 \begin{align*}
 &\hat{f}_n(\bx) \geq 0 \ \text{implies} \ \bar{f}_n(\bx) \geq \frac{1}{2} \ \text{and}\\
 & \hat{f}_n(\bx) < 0 \ \text{implies} \ \bar{f}_n(\bx) < \frac{1}{2}
 \end{align*}
 and consequently we have
 \begin{align*}
  &(1-\bar{f}_n(\bx)) \cdot \mathds{1}_{\{\hat{f}_n(\bx) \geq 0\}}(\bx) + \bar{f}_n(\bx) \cdot \mathds{1}_{\{\hat{f}_n(\bx) <  0\}}(\bx)\\
  = & \min\{1-\bar{f}_n(\bx), \bar{f}_n(\bx)\}\\
  \leq & (1-\bar{f}_n(\bx)) \cdot \mathds{1}_{\{\eta(\bx) \geq \frac{1}{2}\}}(\bx) + \bar{f}_n(\bx) \cdot \mathds{1}_{\{\eta(\bx) <  \frac{1}{2}\}}(\bx). 
 \end{align*}
 Set $g(z) = \log\frac{z}{1-z}$ with $z \in (0,1)$.  Using that 
 \begin{align*}
 h_1(z) = \varphi(1 \cdot g(z)) &= \log\left(1+\exp\left(-\log\frac{z}{1-z}\right)\right)\\
 &= \log\left(1+\frac{1-z}{z}\right)\\
 &= \log\left(\frac{1}{z}\right) = - \log(z)
 \end{align*}
 and 
 \begin{align*}
 h_2(z) = \varphi(-1 \cdot g(z)) &= \log\left(1+\frac{z}{1-z}\right)\\
 &= \log\left(\frac{1}{1-z}\right)\\
 &= - \log(1-z), 
 \end{align*}
 we have
 \begin{align*}
h_1'(z) = - \frac{1}{z} \ \text{and} \ h_2'(z) = \frac{1}{1-z} 
\end{align*}
and consequently
\begin{align*}
|h_1'(z)| = \frac{1}{|z|} \geq 1 \ \text{and} \ |h_2'(z)| = \frac{1}{|1-z|} \geq 1 \ \text{for} \ z \in (0,1). 
\end{align*}
 Using mean value theorem it follows for $j \in \{1, 2\}$ and $z_1, z_2 \in [0,1]$
\begin{align*}
|h_j(z_1) - h_j(z_2)| \geq 1 \cdot |z_1 - z_2|.
\end{align*}
Choosing $z_1=\bar{f}_n(\bx)$ and $z_2 = \eta(\bx)$, this, in turn, means that \eqref{le2eq1} can be further bounded by 
 \begin{align}
 \label{le2eq2}
 2 \cdot \EXP\left\{\left|\varphi\left(Y \cdot g(\bar{f}_n(\bX))\right) - \varphi\left(Y \cdot g(\eta(\bX))\right)\right|\right\}.
 \end{align}
  As $f_{\varphi}^*(\bx) = g(\eta(\bx))$ and $\hat{f}(\bx) = g(\bar{f}_n(\bx))$ by definition and 
 \begin{align*}
 |a-b| \leq |a|+|b|=a+b \quad \mbox{for } a,b \geq 0 
 \end{align*}
 we can then bound \eqref{le2eq2} by
 \begin{align*}
 & 2 \cdot \EXP\left\{\left|\varphi\left(Y \cdot g(\bar{f}_n(\bX))\right) - \varphi\left(Y \cdot g(\eta(\bX))\right)\right|\right\}\\
 & \leq 2  \cdot \EXP\left\{\varphi\left(Y \cdot \hat{f}_n(\bX)\right) + \varphi\left(Y \cdot f_{\varphi}^*(\bX)\right)\right\}\\
 & = 2  \cdot \EXP\left\{\varphi\left(Y \cdot \hat{f}_n(\bX)\right) - \varphi\left(Y \cdot f_{\varphi}^*(\bX)\right)\right\}+ 4 \cdot  \EXP\left\{ \varphi\left(Y \cdot f_{\varphi}^*(\bX)\right)\right\}.
 \end{align*}
    \hfill $\Box$

\subsection{Proofs of Section 3}

This section contains the missing proofs of Theorem 1 and Lemma 2.
\\
\\
  \noindent
        \textbf{Proof of Theorem 1.}
        As 
          \begin{eqnarray*}
          &&
          \EXP\{\mathcal{E}(\hat{f}_n,  f^*)\}
\leq
\sqrt{2} \cdot \EXP\{\mathcal{E}^{\varphi}(\hat{f}_n, f_{\varphi}^*)^{1/2}\}
\leq
\sqrt{2}
\cdot
\EXP\{\mathcal{E}^{\varphi}(\hat{f}_n, f_{\varphi}^*)\}^{1/2}
\end{eqnarray*}
by Lemma 1a),     
     \begin{eqnarray*}
  &&
       \EXP\{\mathcal{E}(\hat{f}_n,  f^*)\}
\leq
2 \cdot \EXP\{\mathcal{E}^{\varphi}(\hat{f}_n, f^*_{\varphi})\}
+ 4 \cdot \mathcal{E}^{\varphi}(f_{\varphi}^*)
  \end{eqnarray*}   
by Lemma 1b) and 
\begin{align}
\label{riskf}
 \mathcal{E}^{\varphi}(f_{\varphi}^*) \leq  \frac{c_{12} \cdot \log n}{\sqrt{n}}.
\end{align}
by using Assumption 2 together with Lemma 3 in \cite{KOK19}, where we choose $F_n =\frac{1}{2} \cdot \log n$,
it is enough to prove that 
        \begin{equation}
          \label{pth1eq1}
          \EXP\{\mathcal{E}^{\varphi}(\hat{f}_n, f_{\varphi}^*)\} 
\leq
c_{13} \cdot (\log n)^{2} \cdot n^{-\min\{\frac{p}{2p+4},
   \frac{1}{4} \}}.
          \end{equation}
Application of Lemma 2, Lemma 8 and Theorem 2
yields
\begin{eqnarray*}
  &&
\EXP\{\mathcal{E}^{\varphi}(\hat{f}_n, f_{\varphi}^*)\}
\\
&&
\leq
2
\cdot
\EXP \left\{
\sup_{f \in \F_{n, \Theta}}
\left|\mathcal{E}^{\varphi}(f) - \mathcal{E}^{\varphi}_n(f)
\right|
\right\}
+
\inf_{f \in \F_{n, \Theta}} \mathcal{E}^{\varphi}(f, f_\varphi^*)
\\
&&
\leq
c_{14} \cdot (\log n)^2 \cdot \frac{\max\{L_n^{(1)},L_n^{(2)}\}}{\sqrt{n}}
+
c_{15} \cdot \left(
\frac{\log L_n^{(2)}}{L_n^{(2)}} + \frac{1}{(L_n^{(1)})^{2p/4}}
\right)
\\
&&
\leq
c_{14} \cdot (\log n)^{2} \cdot \frac{L_n^{(2)}}{\sqrt{n}}
+
c_{15} \cdot 
\frac{\log L_n^{(2)}}{L_n^{(2)}}
+
c_{14} \cdot (\log n)^{2} \cdot \frac{L_n^{(1)}}{\sqrt{n}}
+
c_{15} \cdot \frac{1}{(L_n^{(1)})^{2p/4}}
\\
&&
\leq
c_{16} \cdot (\log n)^{2} \cdot n^{-\min\{\frac{p}{2p+4},
   \frac{1}{4} \}}.
  \end{eqnarray*}

        \hfill $\Box$

\noindent
    \textbf{Proof of Lemma 2.}
This result is the standard error bound for empirical risk
minimization.
For the sake of completeness we present nevertheless a
complete proof.

Let $f \in \F_n$ be arbitrary. Then the definition of $\hat{f}_n$
implies
\begin{align*}
\EXP\{\mathcal{E}^{\varphi}(\hat{f}_n, f^*_{\varphi})\} 
&
\leq
\EXP\{\mathcal{E}^{\varphi}(\hat{f}_n) - \mathcal{E}_n^{\varphi}(\hat{f}_n)
+
\mathcal{E}_n^{\varphi}(\hat{f}_n)-\mathcal{E}_n^{\varphi}(f)
+
\mathcal{E}_n^{\varphi}(f)
-
\mathcal{E}^{\varphi}(f)
+
\mathcal{E}^{\varphi}(f)
-
\mathcal{E}^{\varphi}( f_\varphi^*)\}
\\
&
\leq
\EXP\{\mathcal{E}^{\varphi}(\hat{f}_n) - \mathcal{E}_n^{\varphi}(\hat{f}_n)
+
0
+
\mathcal{E}_n^{\varphi}(f)
-
\mathcal{E}^{\varphi}(f)
+
\mathcal{E}^{\varphi}(f, f_\varphi^*)\}
\\
&
\leq
2
\cdot \EXP\{
\sup_{g \in \F_n}
\left|
\mathcal{E}^{\varphi}(g)
-
\mathcal{E}_n^{\varphi}(g)
\right|\}
+
\mathcal{E}^{\varphi}(f,f_\varphi^*).
\end{align*}

  \quad  \hfill $\Box$
\subsection{Proofs of Section 4}
This section presents the missing proofs of Lemma 3 and Theorem 2.
\\
\\
\noindent
    \textbf{Proof of Lemma 3.}
    For $k \in \{-1,0, \dots, K+1\}$ define
    \[
    B_k(z) =
    \begin{cases}
      0 & ,z < \frac{k-1}{K} \\
     K \cdot( z- \frac{k-1}{K})  & ,\frac{k-1}{K} \leq z < \frac{k}{K} \\
     K \cdot (\frac{k+1}{K}-z)  & ,\frac{k}{K} \leq z < \frac{k+1}{K} \\
     0 &, z \geq \frac{k+1}{K}.
      \end{cases}
    \]
    Note that this implies $B_k(k/K)=1$ and $B_k(j/K)=0$ for $j \in \Z \setminus \{k\}$. Set
    \begin{align*}
    \bar{g}(z)=&g(1/K) \cdot (B_{-1}(z)+B_0(z)) + \sum_{k=1}^{K-1} g(k/K) \cdot B_k(z)+ g(1-1/K) \cdot (B_K(z)+ B_{K+1}(z))\\
    =:&
    \sum_{k=-1}^{K+1} a_k \cdot B_k(z).
    \end{align*}
    Then $\bar{g}$ interpolates the points $(-1/K, g(1/K))$, $(0,g(1/K))$, $(1/K,g(1/K))$,
    $(2/K,g(2/K))$, \dots, 
    $((K-1)/K,g((K-1)/K))$, $(1,g((K-1)/K))$ and $(1+1/K,g((K-1)/K))$,
    is zero outside of $(-2/K,1+2/K)$
    and linear on each interval $[k/K, (k+1)/K]$
    $(k \in \{-2, \dots, K+1\})$.
    Because of
    \[
    B_k(z)=\sigma \left( K \cdot\left( z- \frac{k-1}{K}\right)\right) - 2
    \cdot \sigma \left( K \cdot \left(z-\frac{k}{K}\right)\right) +
    \sigma\left(K \cdot \left(z-\frac{k+1}{K}\right)\right),
    \]
    $\bar{g}$ can be computed by a neural network with
    ReLU activation function and $K+3$ hidden layers with
    $7$ neurons per layer. To see this we use the identity
    \[
    \sigma(x)-\sigma(-x)=\max\{x,0\}-\max\{-x,0\}=x \quad
    \mbox{for } x \in \R,
    \]
    which enables us to compute $\bar{g}$ recursively as follows:
    \[
    \bar{g}(x) = \bar{g}_1^{K+3}-\bar{g}_2^{K+3}
    \]
    where
    \[
    \bar{g}_1^l = \sigma( \bar{g}_1^{l-1} - \bar{g}_2^{l-1})
    + a_{l-2} \cdot B_{l-2}(\bar{g}_3^{l-1} - \bar{g}_4^{l-1})
    \quad
    \mbox{for } l \in \{2, \dots, K+3\},
    \]
    \[
    \bar{g}_2^l = \sigma( - \bar{g}_1^{l-1} + \bar{g}_2^{l-1})
    \quad
    \mbox{for } l \in \{2, \dots, K+3\},
    \]
    \[
    \bar{g}_3^{l} = \sigma( \bar{g}_3^{l-1} - \bar{g}_4^{l-1})
     \quad
    \mbox{for } l \in \{2, \dots, K+3\},
    \]
    \[
\bar{g}_4^{l} = \sigma( - \bar{g}_3^{l-1} + \bar{g}_4^{l-1})
     \quad
    \mbox{for } l \in \{2, \dots, K+3\}
    \]
    and
    \[
    \bar{g}_1^1=a_{-1} \cdot B_{-1}(x), \bar{g}_2^1=0, \bar{g}_3^1=\sigma(x)
    \quad \mbox{and} \quad \bar{g}_4^1=\sigma(-x).
    \]
    By induction it is easy to see that the above recursion implies
    \[
    \bar{g}_3^{l} - \bar{g}_4^{l}=x
    \quad \mbox{and} \quad
     \bar{g}_1^{l}-\bar{g}_2^{l}
     =
     \sum_{k=-1}^{l-2} a_k \cdot B_k(x)
    \]
    for $l \in \{1, \dots, K+3\}$.
    
   Set
\[
h_1(z)=\varphi(g(z)) = \log \left( 1 + \exp \left(
- \log \frac{z}{1-z} \right) \right)
= \log \left( 1 + \frac{1-z}{z}  \right) = - \log z
\]
and
\[
h_2(z)=\varphi(-g(z))
=
 \log \left( 1 + \exp \left(
 \log \frac{z}{1-z} \right) \right)
 =
 - \log (1-z).
\]
In case that  $\eta(\bx) \in [0,2/K]$,  we have
$g(\eta(\bx)) \leq g(2/K)=- \log(K/2-1)<0$.
It further holds that $-1/K \leq \bar{\eta}(\bx) \leq 3/K$ and
$- \log (K-1)=f(1/K) \leq \bar{g}(\bar{\eta}(\bx)) \leq f(3/K) = - \log (K/3-1)$.
Consequently we get 
\begin{align*}
\left|
\eta(\bx) \cdot
\left(
\varphi(\bar{g}(\bar{\eta}(\bx))-
\varphi(g(\eta(\bx))
\right)
\right|
&\leq  \eta(\bx) \cdot \varphi(\bar{g}(\bar{\eta}(\bx)))
+
\eta(\bx) \cdot h_1(\eta(\bx)) \\
& \leq 
\frac{2}{K} \cdot \log(1+\exp(\log(K-1)))
+
\eta(\bx) \cdot \log (\frac{1}{\eta(\bx)}) \\
& \leq  4 \cdot \frac{\log K}{K}, 
\end{align*}
where we use that $z \cdot \log(1/z) \leq (2/K) \cdot \log (K/2)$ for $0<z<2/K$. Accordingly, 
\begin{eqnarray*}
  &&
\left|
(1-\eta(\bx)) \cdot
\left(
\varphi(-\bar{g}(\bar{\eta}(\bx)))-
h_2(\eta(\bx))
\right)
\right|
\\
&&\leq  \varphi(-\bar{g}(\bar{\eta}(\bx)))
+
\varphi(-g(\eta(\bx))) \\
&&=
\log(1+\exp(\bar{g}(\bar{\eta}(\bx))))
+
\log(1+\exp(g(\eta(\bx))))
\\
&&\leq
\log (1 + \exp(- \log(K/3-1)))
+
\log(1+\exp(-\log(K/2-1)))
\\
&& \leq  2 \cdot \exp(-\log(K/3-1)) = \frac{
6
}{
K-3
},
\end{eqnarray*}
where we use in the last inequality that $\log(1+x) \leq x$ for $x >0$.
In a similar way,  it holds that 
\begin{align*}
  \left|
\eta(\bx) \cdot
\left(
\varphi(\bar{g}(\bar{\eta}(\bx))-
\varphi(g(\eta(\bx))
\right)
\right|
+
\left|
(1-\eta(\bx)) \cdot
\left(
\varphi(-\bar{g}(\bar{\eta}(\bx))-
\varphi(-g(\eta(\bx))
\right)
\right|
\leq
10 \cdot \frac{\log K}{K-3}
\end{align*}
for $\eta(\bx) \geq 1-2/K$.
Hence it suffices to show
\begin{align*}
  &
\sup_{\bx \in \Rd, \atop \eta(\bx) \in [2/K,1-2/K]}
\Bigg(
\left|
\eta(\bx) \cdot
\left(
\varphi(\bar{g}(\bar{\eta}(\bx))-
\varphi(g(\eta(\bx))
\right)
\right|
\\
&
\hspace*{7cm}
+
\left|
(1-\eta(\bx)) \cdot
\left(
\varphi(-\bar{g}(\bar{\eta}(\bx))-
\varphi(-g(\eta(\bx))
\right)
\right|
\Bigg)
\\
&
\leq
c_{17} \cdot (\frac{\log K}{K} + \epsilon).
\end{align*}

By the monotonicity of $g$, $|g^\prime(z)|=\frac{1}{z \cdot (1-z)} \geq 1$ for
$z \in (0,1)$, the intermediate value theorem and the definition of $\bar{g}$, 
    we can conclude that for any $\bx \in \Rd$ with
    $\eta(\bx) \in [2/K,1-2/K]$
    we find $\xi_x, \delta_{\bx} \in \R$
    with $|\xi_x| \leq \frac{1}{K}$,
    $|\delta_{\bx}| \leq \frac{1}{K} + \epsilon$ and
    $\eta(\bx)+\delta_{\bx} \in [1/K,1-1/K]$
    such that 
    \begin{equation}
      \label{ple7eq1}
      \bar{g}(\bar{\eta}(\bx))
      =
      g(\bar{\eta}(\bx)+\xi_x)
      =
      g(\eta(\bx) + \delta_{\bx}).
    \end{equation}
This implies
\begin{align*}
  &
\sup_{\bx \in \Rd, \atop \eta(\bx) \in [2/K,1-2/K]}
\Bigg(
\left|
\eta(\bx) \cdot
\left(
\varphi(\bar{g}(\bar{\eta}(\bx))-
\varphi(g(\eta(\bx))
\right)
\right|
\\
&
\hspace*{6cm}
+
\left|
(1-\eta(\bx)) \cdot
\left(
\varphi(-\bar{g}(\bar{\eta}(\bx))-
\varphi(-g(\eta(\bx))
\right)
\right|
\Bigg)
\\
&
\quad
=
\sup_{\bx \in \Rd, \atop \eta(\bx) \in [2/K,1-2/K]}
\Bigg(
|\eta(\bx)| \cdot |h_1(\eta(\bx)+\delta_{\bx})-h_1(\eta(\bx))|
\\
&\hspace*{6cm}
+
|1-\eta(\bx)| \cdot |h_2(\eta(\bx)+\delta_{\bx})-h_2(\eta(\bx))|
\Bigg).
\end{align*}
Consequently it suffices to show that there exists a constant
$c_{18}$ such that we have for any $z \in [2/K,1-2/K]$ and any $\delta \in \R$
with $|\delta| \leq \frac{1}{K} + \epsilon$ and $z+\delta \in [1/K,1-1/K]$
\begin{equation}
  \label{ple7eq2}
  |z| \cdot |h_1(z+\delta)-h_1(z)| \leq c_{18} \cdot
  \left( \frac{1}{K} + \epsilon \right)
  \end{equation}
and
\begin{equation}
  \label{ple7eq3}
    |1-z| \cdot |h_2(z+\delta)-h_2(z)| \leq c_{18} \cdot
  \left( \frac{1}{K} + \epsilon \right).
  \end{equation}
We have
\[
h_1^\prime (z) = - \frac{1}{z}.
\]
By the mean value theorem we get for some $\xi \in [\min\{ z+\delta, z\},
  \max\{ z+\delta,z\}] $
\[
|z| \cdot |h_1(z+\delta)-h_1(z)|
=
|z| \cdot \frac{1}{|\xi|} \cdot |\delta|
\leq
4 \cdot |\delta|
\leq
4 \cdot
  \left( \frac{1}{K} + \epsilon \right),
\]
where we have used that $z, z+\delta \in [1/K,1-1/K]$ and $\delta \leq 2/K$
imply $4 |\xi| \geq |z|$.

In the same way we get
\[
h_2^\prime (z) = \frac{1}{1-z}
\]
and
\[
|1-z| \cdot |h_2(z+\delta)-h_2(z)|
=
|1-z| \cdot \frac{1}{|1-\xi|} \cdot |\delta|
\leq
4 \cdot |\delta|
\leq
4 \cdot
  \left( \frac{1}{K} + \epsilon \right).
\]
\quad \hfill $\Box$

In order to prove Theorem 2 we need the following
three auxiliary results.

\begin{lemma}
  \label{le6}
    Let $d \in \N$,
  let $f:\Rd \rightarrow \R$ be $(p,C)$--smooth for some $p=q+s$,
  $q \in \N_0$  and $s \in (0,1]$, and $C>0$. Let $A \geq 1$
    and $M \in \N$ sufficiently large (independent of the size of $A$, but
     \begin{align*}
       M \geq 2 \ \mbox{and} \ M^{2p} \geq c_{19} \cdot \left(\max\left\{A, \|f\|_{C^q([-A,A]^d)}
       \right\}\right)^{4(q+1)},
     \end{align*}
     where
     \[
     \|f\|_{C^q([-A,A]^d)}
     =
     \max_{\bm{\alpha} \in \N_0^d,  \atop \|\bm{\alpha}\|_1 \leq q}
     \left\|\partial^{\bm{\alpha}} f\right\|_{\infty, [-A,A]^d}
     ,
     \]
 must hold for some sufficiently large constant $c_{19} \geq 1$).
 \\
a) Let $L, r \in \N$ be such that
\begin{enumerate}
\item $L \geq 5+\lceil \log_4(M^{2p})\rceil \cdot \left(\lceil \log_2(\max\{q, d\} + 1\})\rceil+1\right)$
\item $r \geq 2^d \cdot 64 \cdot \binom{d+q}{d} \cdot d^2 \cdot (q+1) \cdot M^d$
\end{enumerate}
hold.
There exists a feedforward neural network
$f_{net, wide}$ with ReLU activation function, $L$ hidden layers
and $r$ neurons per hidden layer such that
\begin{align}
 \| f-f_{net, wide}\|_{\infty, [-A,A]^d} \leq
  c_{20} \cdot \left(\max\left\{A, \|f\|_{C^q([-A,A]^d)}\right\} \right)^{4(q+1)} \cdot M^{-2p}.
  \label{le3eq1}
\end{align}
b) Let $L, r \in \N$ be such that
\begin{enumerate}
\item $L \geq 5M^d+\left\lceil \log_4\left(M^{2p+4 \cdot d \cdot (q+1)} \cdot e^{4 \cdot (q+1) \cdot (M^d-1)}\right)\right\rceil  \lceil \log_2(\max\{q,d\}+1)\rceil+\lceil \log_4(M^{2p})\rceil$
\item $r \geq 132 \cdot 2^d\cdot   \lceil e^d\rceil
  \cdot \binom{d+q}{d} \cdot \max\{ q+1, d^2\}$
\end{enumerate}
hold.
There exists a feedforward neural network
$f_{net, deep}$ with ReLU activation function, $L$ hidden layers
and $r$ neurons per hidden layer such that
 (\ref{le3eq1}) holds with
$f_{net,wide}$ replaced by $f_{net,deep}$.
\end{lemma}

\noindent
    \textbf{Proof.} See Theorem 2 in \cite{KL20}. \hfill $\Box$

\begin{lemma}
  \label{le8}
    Let $d_1,d_2,l \in \N$ with $2^l \leq \min\{d_1,d_2\}$ and set $I=2^l -1 \times 2^l-1$.  
  Define $m$ and $\bar{m}$ by
  \[
m(\bx)=
\max_{
  (i,j) \in \Z^2 \, : \,
  (i,j)+I \subseteq d_1 \times d_2
}
f\left(
\bx_{(i,j)+I}
\right)
\]
and
\[
\bar{m}(\bx)=
\max_{
  (i,j) \in \Z^2 \, : \,
  (i,j)+I \subseteq d_1 \times d_2
}
\bar{f}\left(
\bx_{(i,j)+I}
\right),
\]
where $f$ and $\bar{f}$ satisfy a hierarchical model of level $l$ according to Definition  \ref{de2}. Denote by $g_{k,s}, \bar{g}_{k,s}: \R^4 \to [0,1]$ $(k=1, \dots, l, s=1, \dots, 4^l-k)$ the functions occuring in the definition of $f$ and $\bar{f}$,  respectively.  

Assume that all functions
$g_{k,s}$ $(k=1, \dots, l, s=1, \dots, 4^l-k)$ are Lipschitz continuous regarding the Euclidean distance with Lipschitz constant $C>0$ and that
\begin{equation}
\|\bar{g}_{k,s}\|_{[-2,2]^4,\infty}\leq2.
\label{le4eq3}
\end{equation}
Then
\[
|m(\bx)-\bar{m}(\bx)| \leq
(2C+1)^l
\cdot \max_{i \in \{1, \dots,l\}, s \in \{ 1, \dots, 4^{l-i}\}}
\|                      g_{i,s}
                -
                        \bar{g}_{i,s}
\|_{[-2,2]^4,\infty}
\]
for any
$\bx \in [0,1]^{d_1 \times d_2}$.
  \end{lemma}

\noindent
    \textbf{Proof.} The assertion follows from Lemma 4 in \cite{KKW22}. \hfill $\Box$

\begin{lemma}
  \label{le9}
   Let $d_1,d_2,l\in \N$ with $2^l\leq\min\{d_1,d_2\}$. For $k \in \{1, \dots,l\}$
  and $s \in \{1, \dots, 4^{l-k}\}$ let
  \[
\bar{g}_{net, k,s} : \R^4 \rightarrow \R
\]
be defined by a feedforward neural network with $L_{net}\in\N$
hidden layers and $r_{net}\in\N$ neurons per hidden layer and ReLU
activation function.
Set
\[I=2^{l}-1 \times 2^{l}-1\]
 and
define
$\bar{m}:[0,1]^{d_1 \times d_2} \rightarrow \R$
by
\[
\bar{m}(\bx)=
\max_{
  (i,j) \in \Z^2 \, : \,
  (i,j)+I \subseteq d_1 \times d_2
}
\bar{f}\left(
\bx_{(i,j)+I}
\right),
\]
where $\bar{f}$ satisfies
    \[
\bar{f}=\bar{f}_{l,1}
    \]
    for some
    $\bar{f}_{k,s} :[0,1]^{2^k  \times 2^k} \rightarrow \R$ recursively defined by
    \begin{eqnarray*}
    \bar{f}_{k,s}(\bx)&=&\bar{g}_{net,k,s} \big(
    \bar{f}_{k-1,4 \cdot (s-1)+1}(\bx_{
1,1})
    , 
        \bar{f}_{k-1,4 \cdot (s-1)+2}(\bx_{2,1}), 
        \bar{f}_{k-1,4 \cdot (s-1)+3}(\bx_{
1,2
        }), 
        \bar{f}_{k-1,4 \cdot s}(\bx_{
2,2
        })
    \big)
    \end{eqnarray*}
    for $k=2, \dots, l, s=1, \dots,4^{l-k}$, $\mathbf{x} \in [0,1]^{2^k \times 2^k}$
 and
    \[
 \bar{f}_{1,s}(x_{1,1},x_{1,2},x_{2,1},x_{2,2})= \bar{g}_{net,1,s}(x_{1,1},x_{1,2},x_{2,1},x_{2,2})
 \]
 for $s=1, \dots, 4^{l-1}$.
 Set
 \[
l_{net}=\frac{4^{l}-1}{3} \cdot L_{net}+l,
\]
\[
k_s=\frac{2 \cdot 4^{l} + 4}{3} +r_{net}
\quad
(s=1, \dots, l_{net}),
\]
and set
\[M_s=2^{\pi(s)}\quad \mbox{for }s\in\{1,\dots,l_{net}\},\]
where the function $\pi:\{1,\dots,l_{net}\}\rightarrow\{1,\dots,l\}$ is defined by 
\[\pi(s)=\sum_{i=1}^{l}\IND_{\left\{s\geq i+\sum_{r=l-i+1}^{l-1}4^r\cdot L_{net}\right\}}.\]
 Then there exists some $m_{net} \in \F^{(CNN)}_{l_{net},\bk,\bM}$ such that
 \[
\bar{m}(\bx) = m_{net}(\bx)
 \]
holds for all $\bx \in [0,1]^{\{1, \dots, d_1\} \times \{1, \dots, d_2\}}$.
  \end{lemma}

\noindent
    \textbf{Proof.} See Lemma 5 in \cite{KKW22}. \hfill $\Box$

\noindent
    \textbf{Proof of Theorem 2.}
    For each $g_{k,s}$ in the hierarchical max--pooling model
    for $\eta$ we select an approximating neural network from
    Lemma \ref{le6} a) which approximates $g_{k,s}: \R^4 \rightarrow \R$
    up to an error
    of order $(L_n^{(1)})^{-2p/4}$. Then we use Lemma \ref{le9} to generate
    with these networks a convolutional neural network, and
    combine this network with the feedforward neural network with $L_n^{(2)}$
    layers and $7$ neurons per layer of Lemma 4.
    For the corresponding network $h \in \F_n$ we get by Lemma 4
    and Lemma \ref{le8}
\begin{align*}
&\sup_{x \in \R^{d_1 \times d_2}}
\left(
\left|
\eta(\bx) \cdot
\left(
\varphi(h(\bx))-
\varphi(f_\varphi^*(\bx))
\right)
\right|
+
\left|
(1-\eta(\bx)) \cdot
\left(
\varphi(-h(\bx))-
\varphi(-f_\varphi^*(\bx))
\right)
\right|
\right)
\\
&
\leq
c_{21} \cdot \left(\frac{\log L_n^{(2)}}{L_n^{(2)}} + \frac{1}{(L_n^{(1)})^{2p/4}}\right).
\end{align*}
Because of
\begin{align*}
\inf_{f \in \F_{n, \Theta}}
\mathcal{E}^{\varphi}(f, f_{\varphi}^*)
&
\leq
\mathcal{E}^{\varphi}(h, f_{\varphi}^*)
\\
&
=
\int
(
\eta(\bx) \cdot \varphi( h(\bx))
+
(1-\eta(\bx)) \cdot \varphi(- h(\bx))
)
\PROB_\bX(d\bx)
\\
&
\quad
-
\int
(
\eta(\bx) \cdot \varphi( f_\varphi^*(\bx))
+
(1-\eta(\bx)) \cdot \varphi(- f_\varphi^*(\bx))
)
\PROB_\bX(d\bx)
\\
&
\leq
\sup_{\bx \in [0,1]^{d_1 \times d_2}}
\Big(
\left|
\eta(\bx) \cdot
\left(
\varphi(h(\bx))-
\varphi(f_\varphi^*(\bx))
\right)
\right|
\\
&
\quad
+
\left|
(1-\eta(\bx)) \cdot
\left(
\varphi(-h(\bx))-
\varphi(-f_\varphi^*(\bx))
\right)
\right|
\Big)
\end{align*}
this implies the assertion.
    
    \hfill $\Box$

\subsection{Proofs of Section 5}
This section contains the bound and corresponding proof of the generalisation error $\sup_{f \in \mathcal{F}_{n, \Theta}} |\mathcal{E}^{\varphi}(f) - \mathcal{E}_n^{\varphi}(f)|$.

\begin{lemma}
  \label{le3}
  Let $\varphi$ be the logistic loss and $\mathcal{D}_n= \{(\bX_i, Y_i)\}_{i=1}^n$ . Let $\F_{n, \Theta}$ be defined as in (5). Then
\begin{eqnarray*}
&&\EXP \left\{
\sup_{f \in \mathcal{F}_{n, \Theta}} |\mathcal{E}^{\varphi}(f) - \mathcal{E}_n^{\varphi}(f)|
\right\} \leq  c_{22} \cdot (\log n)^2 \cdot \frac{\max\{L_n^{(1)},L_n^{(2)}\}}{\sqrt{n}}.
\end{eqnarray*}
  \end{lemma}

\noindent
    \textbf{Proof of Lemma \ref{le3}.}
Since $f \in \mathcal{F}_{n, \Theta}$ satisfies $\|f\|_\infty \leq \beta_n = c_1 \cdot \log n$, we have for any $\bx \in \R^{d_1 \times d_2}$ and $y \in \{-1,1\}$ 
\begin{align*}
\varphi(y \cdot f(\bx)) =\log(1+e^{-y \cdot f(\bx)}) \leq \log(1+e^{|f(\bx)|}) \leq \log(1+e^{\beta_n}) = c_{23} \cdot \log n.
\end{align*}
Set 
\begin{align*}
\mathbf{Z}=(\bX, Y), \mathbf{Z}_1=(\bX_1, Y_1), \dots, \mathbf{Z}_n=(\bX_n, Y_n), 
\end{align*}
and 
\begin{align*}
  \mathcal{H}_n= \{h: \R^{d_1 \times d_2} \times \R \to \R: \exists f \in
  \mathcal{F}_n \ \text{such that} \ h(\bx,y) = \varphi(y\cdot f(\bx))\}.
\end{align*}
 By Theorem 9.1 in \cite{GKKW02} one has, for arbitrary $\epsilon > 0$   
 \begin{align*}
\PROB\left\{\sup_{f \in \mathcal{F}_{n, \Theta}} |\mathcal{E}^{\varphi}(f) - \mathcal{E}_n^{\varphi}(f)| > \epsilon\right\}
& = \PROB\left\{\sup_{h \in \mathcal{H}_n} \left|\EXP h(\mathbf{Z}) - \frac{1}{n} \sum_{i=1}^n h(\mathbf{Z}_i)\right| > \epsilon\right\}\\
& \leq 8 \EXP\left\{\mathcal{N}_1\left(\frac{\epsilon}{8}, \mathcal{H}_n, \mathbf{Z}_1^n\right)\right\}e^{-\frac{n \epsilon^2}{128 \cdot c_{23}^2 \cdot (\log n)^2}}.
 \end{align*}
 Let
 $h_i(\bx, y) = \varphi(y \cdot f_i(\bx))$ $((\bx,y) \in \R^{d_1 \times d_2} \times \{-1,1\})$ for some
 $f_i: \R^{d_1 \times d_2} \rightarrow \R$.  Together with the Lipschitz continuity of $\varphi$,  it follows
\begin{align*}
\frac{1}{n} \sum_{i=1}^n \left|h_1(\mathbf{Z}_i) - h_2(\mathbf{Z}_i)\right|
&= \frac{1}{n} \sum_{i=1}^n \left|\varphi(Y_i \cdot f_1(\bX_i)) - \varphi(Y_i \cdot f_2(\bX_i))\right|\\
& \leq \frac{1}{n} \sum_{i=1}^n \left|Y_i f_1(\bX_i) - Y_i f_2(\bX_i)\right|\\
& = \frac{1}{n} \sum_{i=1}^n \left|f_1(\bX_i) - f_2(\bX_i)\right|.
\end{align*}   
Thus, if $\{f_1, \dots, f_{\ell}\}$ is a $\epsilon$-cover of  $\mathcal{F}_n$, then $\{h_1, \dots, h_{\ell}\}$ is a $\epsilon$-cover of $\mathcal{H}_n$. Then
\begin{align*}
\mathcal{N}_1\left(\frac{\epsilon}{8}, \mathcal{H}_n, \mathbf{Z}_1^n\right) \leq \mathcal{N}_1\left(\frac{\epsilon}{8}, \mathcal{F}_n, \mathbf{X}_1^n\right).
\end{align*}
Set $L_{max}=\max\{L_n^{(1)},L_n^{(2)}\}$.
Lemma 4 implies
 \begin{align*}
   \mathcal{N}_1\left(\frac{\epsilon}{8}, \mathcal{F}_n, X_1^n\right) \leq
   \left(\frac{c_{24} \cdot \log n}{\epsilon}\right)^{c_{25} \cdot L_{\max}^2 \cdot \log(L_{\max} \cdot d_1 \cdot d_2)}
 \end{align*}
 and
 \begin{align*}
 \label{le3eq1}
 &\PROB\left\{\sup_{f \in \mathcal{F}_n} \left|\mathcal{E}^{\varphi}(f) - \mathcal{E}_n(f)\right| > \epsilon\right\} \leq 8 \cdot \left\{ \left(\frac{ c_{24} \cdot \log n}{\epsilon}\right)^{c_{25} \cdot L_{\max}^2 \cdot \log(L_{\max} \cdot d_1 \cdot d_2)}\right\} e^{-\frac{n \epsilon^2}{c_{23} (\log n)^2}}.
 \end{align*}
 Set
 \begin{eqnarray*}
 \epsilon_n&=& \sqrt{
   2 \cdot
   \frac{c_{23} \cdot (\log n)^2}{n}
   \cdot \log \left(
8 \left\{ \left(\frac{ c_{24} \cdot \log n}{1/n}\right)^{c_{25} \cdot L_{\max}^2 \cdot \log(L_{\max} \cdot d_1 \cdot d_2)}\right\}
   \right)
 }
 \\
 &=&
    c_{26} \cdot (\log n)^{2} \cdot 
      \frac{\max\{L_n^{(1)},L_n^{(2)}\}}{\sqrt{n}}
 \geq \frac{1}{n}.
 \end{eqnarray*}
 Here we have used for the last inequality that w.l.o.g. we can
 assume $\min\{c_{10},c_{12},c_{13}\} \geq 1$.
 Then
    \begin{align*}
&    \EXP \left\{
\sup_{f \in \F_{n, \Theta}}
\left|
\mathcal{E}^{\varphi}(f) - \mathcal{E}_n^{\varphi}(f)
\right|
\right\}\\
&\leq \epsilon_n+ \int_{\epsilon_n}^{\infty}8 \left\{\left(\frac{c_{24} \cdot \log n}{t}\right)^{c_{25} \cdot L_{\max}^2 \cdot \log(L_{\max} \cdot d_1 \cdot d_2)}\right\} e^{-\frac{n t^2}{c_{23} \cdot (\log n)^2}}dt\\
&= \epsilon_n+ \int_{\epsilon_n}^{\infty}8 \left\{\left(\frac{c_{24} \cdot \log n}{t}\right)^{c_{25} \cdot L_{\max}^2 \cdot \log(L_{\max} \cdot d_1 \cdot d_2)}\right\} e^{-\frac{n t^2}{2c_{23} \cdot (\log n)^2}}e^{-\frac{n t^2}{2c_{23} \cdot (\log n)^2}}dt\\
& \leq \epsilon_n +
\int_{\epsilon_n}^{\infty}
e^{-\frac{n t^2}{2 \cdot c_{23} \cdot (\log n)^2}} dt
  \\
  & \leq
  \epsilon_n +
  \int_{\epsilon_n}^{\infty}
  e^{-\frac{n \cdot \epsilon_n \cdot t}{2 \cdot c_{23} \cdot (\log n)^2}} dt
    \\
    & \leq
    c_{27} \cdot (\log n)^{2} \cdot 
      \frac{\max\{L_n^{(1)},L_n^{(2)}\}}{\sqrt{n}}.
    \end{align*}
    
\quad    \hfill $\Box$

\end{document}